\documentclass{amsart}

\usepackage{amssymb}
\usepackage{amsthm}
\usepackage{amsmath}
\usepackage{amsbsy}
   \usepackage{bbm}
    \usepackage[all]{xy}
    \usepackage{bm}
    \usepackage{hyperref}
    \usepackage{placeins}
\usepackage{tikz}
\usepackage{graphicx}
\usepackage{subcaption}
\usepackage{mathtools,xparse}
\usepackage{mathrsfs}
\usepackage{algorithm}
\usepackage{algorithmic}
 \usepackage{cleveref}
\crefname{figure}{Figure}{Figures}
\crefname{algorithm}{Algorithm}{Algorithms}
\crefname{appendix}{Appendix}{Appendices}
\def\deltat{{\Delta t}}
\def\dd{{\rm d}}
\def\x{{\bm x}}

\def\p{{\bm p}}
\def\W{{\bm W}}
\def\0{{\bm 0}}
\def\I{{\bm I}}
\def\numx{\bm x}
\def\nump{\bm p}
\def\Rand{\mathcal{R}}

\DeclarePairedDelimiter{\abs}{\lvert}{\rvert}
\DeclarePairedDelimiter{\norm}{\lVert}{\rVert}
\DeclareMathOperator{\vol}{vol}

\begin{document}

\title{Quadrature Points via Heat Kernel Repulsion}  

\thanks{This work was partially supported by the National Science Foundation under Grant DMS-1638521 to the Statistical and Applied Mathematical Sciences Institute. The research of J.L. was also supported in part by the National Science Foundation under award DMS-1454939.}

\author[]{Jianfeng Lu}
\address[Jianfeng Lu]{Department of Mathematics, Department of Physics, and Department of Chemistry,
Duke University, Box 90320, Durham NC 27708, USA}
\email{jianfeng@math.duke.edu}

\author[]{Matthias Sachs}
\address[Matthias Sachs]{Department of Mathematics, 
Duke University, Box 90320, Durham NC 27708, USA and 
Statistical and Applied Mathematical Sciences Institute (SAMSI), Durham NC 27709, USA}
\email{msachs@math.duke.edu}

\author[]{Stefan Steinerberger}
\address[Stefan Steinerberger]{Department of Mathematics, Yale University, New Haven, CT 06510, USA}
\email{stefan.steinerberger@yale.edu}

\begin{abstract} 
We discuss the classical problem of how to pick $N$ weighted points on a $d-$dimensional manifold so as
to obtain a reasonable quadrature rule
$$ \frac{1}{|M|}\int_{M}{f(x) dx} \simeq \sum_{n=1}^{N}{a_i f(x_i)}.$$
This problem, naturally, has a long history; the purpose of our paper is to propose selecting points and weights
so as to minimize the energy functional
$$ \sum_{i,j =1}^{N}{ a_i a_j \exp\left(-\frac{d(x_i,x_j)^2}{4t}\right) } \rightarrow \min, \quad \mbox{where}~t \sim N^{-2/d},$$
$d(x,y)$ is the geodesic distance and $d$ is the dimension of the manifold. This
yields point sets that are theoretically guaranteed, via spectral theoretic properties of the Laplacian $-\Delta$, to have good properties.  
One nice aspect is that the energy functional is universal and independent of the underlying manifold; we show several numerical examples.
\end{abstract}
\maketitle

\section{Introduction}

\textbf{Introduction.} We study the following problem. Given a compact, connected, and smooth Riemannian manifold $(M,g)$ how would one distribute $N$ weighted points
on $M$ to achieve a good quadrature rule? Or, put differently, if our goal is
$$ \frac{1}{|M|}\int_{M}{f(x) dx} \simeq  \sum_{n=1}^{N}{a_i f(x_i)},$$
how should one distribute the pairs $(a_i, x_i) \in \mathbb{R} \times M$? This question is, of course, classical and entire branches of numerical analysis have
evolved around it (\cite{davis, dick, drmota, kuipers, novak1, novak2, novak3}). The richness of the question stems from its ambiguity since it is not a priori clear in which function class one should assume $f$ to be nor 
how the integration error should be quantified. Two particular ways of approaching the problem (among many others that we do not discuss) are as follows:
\begin{enumerate}
\item pick a set of functions $f_1, f_2, \dots, f_k$ and choose the points in such a way that these functions are integrated exactly 
\item model points as charged particles that repel each other and study minimal energy configurations of the dynamical system.
\end{enumerate}
The hope being,
for the first scenario, that by picking suitable `representative'
functions, one is guaranteed to integrate elements in
$\mbox{span}\left\{f_1, \dots, f_k\right\}$ exactly and would hope
that similar functions should not behave terribly. In light of this,
it makes sense to choose the functions $f_1, \dots, f_k$ to be
orthogonal in $L^2$ and to choose them as smoothly as possible. When
working on $\mathbb{S}^{d-1}$ and using polynomials as test functions,
this gives rise to the notion of spherical design which is an entire
subject in itself (see \cite{bond, bond2, conway, delsarte, grab,
  leb1, leb5, sey, sobolev, yudin}). The second approach dates back to
a 1904 paper of J. J. Thomson where the problem was raised for $N$
equally charged electrons on $\mathbb{S}^2$ (though not for the
purpose of numerical integration but because of intrinsic interest). A
definite classification of minimizers seems out of reach, the case
$N=5$ was only settled very recently \cite{schwartz}. Needless to say,
the second question, even restricted to the sphere $\mathbb{S}^{d-1}$
has given rise to entire subfields too, we refer to the recent survey
of Brauchart \& Grabner \cite{brau} for an overview. The second
question is also of relevance in mathematical physics (see
e.g. \cite{abri, chatterjee, dahlberg, serf}, we refer to a survey of
Blanc \& Lewin on the Crystallization Conjecture \cite{blanc}).

\medskip 

\textbf{Related results.} Our approach here is inspired by the recent paper by one of the authors
\cite{stein}. We start by describing \cite{stein} which combines the
two approaches mentioned in the introduction: clearly, on any given
manifold, a natural set of low-frequency orthogonal functions is given
by the eigenfunctions of the Laplacian $-\Delta_g$. Indeed, by the
Courant-Fischer-Weyl minimax principle, these can be said to be
uniquely characterized by orthogonality and the requirement of
minimizing the Dirichlet functional. It is thus not unreasonable to
measure the quality of a set of weighted points by how well they
integrate the first few Laplacian eigenfunctions.  On the torus, the
Laplacian eigenfunctions are given by
$e^{2\pi i \left\langle k, x \right\rangle}$, where
$k \in \mathbb{Z}^d$.  The question of finding lower bounds for the
expression
$$ \sum_{\|k\|_{\infty} \leq X}{ \left| \sum_{n=1}^{N}{ e^{2\pi i  \left\langle k, x_n \right\rangle}} \right|^2}$$ 
over all possible sets of $N$ points $\left\{x_1, \dots, x_N \right\} \subset \mathbb{T}^d$ arose in the study of irregularities of distribution. A now 
 classical result of Montgomery \cite{mont1, mont} implies that a set of $N$ points cannot be orthogonal to more than the first $\sim c_d N$ trigonometric functions.  Bilyk \& Dai \cite{bilyk} have established the analogous result on $\mathbb{S}^{d-1}$ with trigonometric functions replaced by spherical harmonics (also in the context of irregularities of distribution).
We refer to \cite{stein1} for a recent refinement of Montgomery's result and to recent work of Bilyk, Dai and the third author \cite{neu} for general refinements. The paper \cite{stein} establishes that a set of $N$ points on a compact $d-$dimensional manifold cannot integrate more than the
first $c_d N + o(N)$ eigenfunctions exactly (see also \cite{ahrens, steingraph}). Both papers \cite{neu, stein} make  use of the heat kernel. 

\medskip

\textbf{Applications.} We believe that this points as well as the main idea behind their construction are likely to have applications
in a variety of fields. An obvious example is that of Numerical Integration: the idea of enforcing as much orthogonality to low-frequency
eigenfunctions as possible is well in line with the underlying idea of spherical harmonics. However, a big advantage of our approach
is that by using the heat kernel (or, the Gaussian kernel), we do not need as strong assumptions on the geometry. In particular, our
approach may be very well suited on finite Graphs and, as such, our points may be suitable in several modern problems in data science.
We give one representative: diffusion maps \cite{diff1,diff2} are a popular dimensionality reduction method that constructs a nonlinear
embedding of a manifold using Laplacian eigenfunctions (or Laplacian eigenvectors of the Graph Laplacian in the discrete case); we
refer to Jones, Maggioni \& Schul \cite{p1, p2} for the rigorous theory.
If one were to subsample the manifold, our points have the nice byproduct of preserving an overall global structure. This may also
have implications for the numerical integration of smooth functions on graphs (where smoothness is defined with respect to a Graph Laplacian, see \cite{linderman}).

\medskip 

\textbf{Organization.} The main idea of our proposed construction of the quadrature points
and weights is to minimize an interaction energy motivated from the
above consideration of heat kernels, as will be presented in next
Section. The construction is validated by numerical experiments to
compare the proposed approach with several existing approaches in
Section~\ref{sec:numerics}.

\section{Quadrature points construction}
\subsection{Heat kernel}\label{sec:heat:kernel} Before stating the main idea, we quickly
recall the heat kernel. On a compact manifold as specified above, the Laplace operator
$-\Delta_g$ has discrete spectrum, eigenvalues
$0 = \lambda_0 < \lambda_1 \leq \lambda_2 \leq \dots$ and
eigenfunctions $\phi_0, \phi_1, \dots \in L^2(M)$. These
eigenfunctions diagonalize the heat flow which allows us to explicitly
solve the heat equation for arbitrary initial data $f \in L^2(M)$. It
suffices to note that
\begin{equation}\label{eq:heat:flow}
 u(t,x) = \sum_{n=0}^{\infty}{ e^{-\lambda_n t} \left\langle f, \phi_n\right\rangle \phi_n} \qquad \mbox{satisfies} \quad \partial_t u = 
\Delta_g u
\end{equation}
and, by completeness, $u(0,x) = f(x)$.  The heat kernel
$p(t,x,y):(0,\infty) \times M \times M \rightarrow \mathbb{R}$ is then
given by
$$ p(t,x,\cdot) =\sum_{n=0}^{\infty}{ e^{-\lambda_n t} \left\langle \delta_x, \phi_n\right\rangle \phi_n} = \sum_{n=0}^{\infty}{ e^{-\lambda_n t}  \phi_n(x) \phi_n} .$$
While a precise computation requires the precise solution of the heat equation (or, equivalently, complete knowledge of the Laplacian eigenfunctions), the behavior
for small $t$ is easier to understand (we refer to Hsu \cite{hsu} for details). Varadhan's short time asymptotic \cite{varadhan} states
$$ \left[e^{t\Delta}\delta_x\right](y) =  \frac{1+o(1)}{(4\pi t)^{d/2}}\exp\left(-\frac{d(x,y)^2}{4t}\right) \qquad \mbox{as}~t \rightarrow 0,$$
where $d(x,y)$ is the geodesic distance on the manifold, and $(e^{t\Delta})_{t\geq 0}$ denotes the semigroup associated with the heat equation, i.e., $e^{t\Delta} f(x) = u(t,x)$ with $u$ as defined in \eqref{eq:heat:flow}.

\subsection{Quadrature points.} We now present the main idea: there is a fairly canonical way of finding sets of points that arise as the optimal configuration
of an energy functional while simultaneously being good at integrating low-frequency functions: that functional is given by
\begin{equation*}\label{eq:min:heat:kernel}
\sum_{i,j =1}^{N}{p(t,x_i, x_j)} \rightarrow \min
\end{equation*}
for some $t > 0$ to be determined.  This can be easily seen via an
expansion into Laplacian eigenfunctions and self-adjointness of the
heat propagator
\begin{equation}\label{eq:identity:1}
\begin{aligned}
 \sum_{i,j =1}^{N}{p(t,x_i, x_j)} &= \left\langle \sum_{i=1}^{N}{ e^{t \Delta/2} \delta_{x_i}} , \sum_{i=1}^{N}{e^{t \Delta/2} \delta_{x_i}} \right\rangle = \left\| e^{t\Delta/2}   \sum_{i=1}^{N}{ \delta_{x_i}} \right\|_{L^2}^2 \\
&= \sum_{k=0}^{\infty}{ \left| \left\langle e^{t\Delta/2}   \sum_{i=1}^{N}{ \delta_{x_i}} , \phi_k \right\rangle  \right|^2}= \sum_{k=0}^{\infty}{ \left| \left\langle   \sum_{i=1}^{N}{ \delta_{x_i}} , e^{t\Delta/2}  \phi_k \right\rangle  \right|^2}\\
&= \sum_{k=0}^{\infty}{e^{-\lambda_k t} \left| \left\langle   \sum_{i=1}^{N}{ \delta_{x_i}} ,   \phi_k \right\rangle  \right|^2} =   \sum_{k=0}^{\infty}{e^{-\lambda_k t} 
\left(    \sum_{i=1}^{N}{ \phi_k(x_i)}   \right)^2      } 
\end{aligned}
\end{equation}
Minimizing the heat kernel functional corresponds to minimizing a weighted integration error against the Laplacian eigenfunctions (since all but
$\phi_0$ have mean value 0). This also suggests a natural choice for the value of $t$: one could not, in general, hope to be able to integrate more than the
first $N$ eigenfunctions exacty, which suggests to put $t \sim \lambda_N^{-1}$. Using Weyl's law $\lambda_N \sim N^{2/d}$ then suggests the scaling $t \sim N^{-2/d}$ on a $d-$dimensional
manifold. This is also the natural scaling for nearest neighbor interactions. The main downside of the method is that the computation of $p(t,x,y)$ is generally quite difficult and only known in closed form for a small number
of special manifolds. However, this is where the short-term asymptotic of $p(t,x,y)$ can be favorably employed since $t \sim N^{-2/d} \ll 1$. Altogether, we are thus proposing to minimize (ignoring a multiplicative factor depending on $t$)
\begin{equation}\label{eq:minimisation:energy}
E_{\rm Gaussian}(x) =  \sum_{i,j =1}^{N}{  \exp\left(-\frac{d(x_i,x_j)^2}{4t}\right) } \rightarrow \min \qquad \mbox{where}~t \sim N^{-\frac{2}{d}}.
\end{equation}
The very same method also works for weighted points $(a_i, x_i) \in \mathbb{R} \times M$ and suggests, by the very same reasoning, the minimization of
\begin{equation*}
\widetilde{E}_{\rm Gaussian}(a,x) = \sum_{i,j =1}^{N}{ a_i a_j \exp\left(-\frac{d(x_i,x_j)^2}{4t}\right) } \rightarrow \min  \qquad \mbox{where}~t \sim N^{-\frac{2}{d}}.
\end{equation*}
We emphasize that this approach is quite different from the more classical and popular methods that employ kernels of the type $\|x_i-x_j\|^{-s}$ for some $s > 0$ (Riesz kernels) which yield point configurations minimizing the so-called Riesz energy. Another advantage is that our functional automatically selects suitable weights as a byproduct (though this could conceivably also be implemented for Riesz energies or, generically, for any given set of points).\\

The problem of constructing $N$ support points which minimize the heat kernel functional \ref{sec:heat:kernel}  or its' approximation \eqref{eq:minimisation:energy} can be understood in a physical context as the problem of optimally placing $N$ heat sources on the manifold $M$ such that the difference between the temperature distribution at time $t/2$  and a uniform distribution of temperature is minimized in $L^{2}(M)$. 


\section{Numerical Results}\label{sec:numerics}

\subsection{On the torus.}  The $d-$dimensional torus $\mathbb{T}^d$ vastly exceeds any other manifold in simplicity: in terms of numerical integration, many sets of points have been proposed and
studied. We will compare our approach to both deterministic (1-4) and stochastic (5-6) (Quasi-)Monte Carlo sampling techniques
\begin{enumerate}
\item Minimal energy configurations for Riesz energies (see \cite{brau})
\item Halton sequence \cite{halton}
\item Good Lattice Point constructions \cite{hlawka}
\item Sobol sequence \cite{sobol}
\item i.i.d. uniform random samples
\item Latin Hypercube samples (LHS)
\end{enumerate}

\subsubsection{Construction of quadrature points and weights}\label{sec:torus:construction}
We construct unweighted quadrature points as (approximate) solutions of the non-convex minimisation problem \eqref{eq:minimisation:energy} by employing a simulated annealing scheme based on an underdamped Langevin equation. For the purpose of comparison, approximate minimal energy configurations for Riesz energies
\begin{equation*}
E_{{\rm Riesz}, s}(x) = \sum_{i,j=1}^{N}\frac{1}{d(x_{i},x_{j})^{s}}, s>0,
\end{equation*}
are obtained by the same simulated annealing scheme. We provide details on the annealing scheme in \cref{ap:sec:annealing:torus}.
Weights $\left\{a_1, \dots, a_N\right\}$ for a given set of points $\{x_{1},\dots,x_{N}\}$ are then constructed by solving 
the quadratic minimisation problem
\begin{equation}\label{eq:weight:opt}
a = \arg\min_{a^{\prime} \in \mathbb{R}^{N}} \sum_{i,j =1}^{N}{ a^{\prime} _i a^{\prime} _j \exp\left(-\frac{d(x_i,x_j)^2}{4t}\right) } \quad \mbox{subject to} \quad \sum_{n=1}^{N}{a_n} = \vol(M).
\end{equation}
It is easy to see that the solution of \eqref{eq:weight:opt} takes the closed form
\begin{equation}\label{eq:sol:weights}
a=\frac{ C^{-1} \mathbbm{1}_{N}}{\mathbbm{1}_{N}^{\top}C^{-1} \mathbbm{1}_{N} },
\end{equation}
where $ \mathbbm{1}_{N} = (1,1,\dots,1)^{\top}\in \mathbb{R}^{N}$, 
and 
\[
C := \left [ \exp\left(-\frac{d(x_i,x_j)^2}{4t}\right) \right ]_{1 \leq i,j \leq N} \in \mathbb{R}^{N\times N}.
\]
We note that for all point configurations considered in our numerical experiments the expression  \eqref{eq:sol:weights}  resulted in positive weights $(a_{i})_{1\leq i \leq N}$. In particular, $a/\vol(M)$ is contained in the standard $(N-1)$-simplex. This may be explained by the fact that for sufficiently small off-diagonal terms of $C$, the matrix $C^{-1}$ is diagonally dominant and thus the expression \eqref{eq:sol:weights} is guaranteed to be positive. The computation of points relies on the geodesic distance: since predominant interactions are on the scale $N^{-1/d} \ll 1$ and the objects under consideration are either flat (the torus) or very symmetric (the sphere), we simplified computation by using the Euclidean distance of an embedding. Theoretically, one could obtain slightly better results on the sphere by incorporating the term correcting for curvature but this effect is of lower order.  \\ 
\subsubsection{Evaluation of quadrature error}
As above, let $ 0=\lambda_{0} <  \lambda_{1} \leq \lambda_{2} \dots  $ be the eigenvalues of the Laplace operator $-\Delta$, considered as an operator on $\mathcal{C}^{\infty}(\mathbb{T}^{d},\mathbb{C})$. For any $s \in \mathbb{N}$ we use $\phi_{s}(x) = e^{i k_{s}x}, ~k_{s}\in \mathbb{N}^{d},$ to denote the eigenfunction associated to the eigenvector $\lambda_{s}$ with eigenvalue $\lambda_{s} = \norm{ k_{s}}^{2}$.
For a prescribed (weighted) point set $\left ( a_{j},x_{j} \right )_{1 \leq j \leq N}$, and eigenvalue $\lambda_{s} \in \sigma(-\Delta)$ we abbreviate the integration error of the associated quadrature scheme by
\[
\begin{aligned}
\mathcal{E}_{\lambda_{s}} &:=  \abs*{\sum_{j=1}^{N}a_{j} e^{i k_{s} \cdot x_{j}} }^{2}.
\end{aligned}
\]
Since we are interested not only
in the integration error with respect to single eigenfunctions but
also in the qualitative behaviour of the integration error for
low-frequency eigenfunctions, we will also work with the sum over all
eigenfunctions up to a certain index
\[
\mathcal{E}_{\leq s} := \sum_{l =1}^{s} \mathcal{E}_{\lambda_{l}}.
\]
The Cauchy-Schwarz inequality immediately implies that this is the maximal integration error in the space $V$ spanned by the linear combination of the first $s$ eigenfunctions with respect to that particular quadrature rule. Using the projection $\pi:L^2(\mathbb{T}^d) \rightarrow V$, we see that the normalization of the weights implies that
\begin{equation}\label{eq:inequality:1}
\begin{aligned}
\sup_{f \in V \atop \|f\|_{L^2} \leq 1}\left| \int_{\mathbb{T}^d}{f(x) dx} - \sum_{i=1}^{N}{a_i f(x_i)}\right|^2 &= \sup_{f \in V \atop \|f\|_{L^2} \leq 1}\left| \sum_{i=1}^{N}{a_i f(x_i)}\right|^2 \\
&= \sup_{f \in V \atop \|f\|_{L^2} \leq 1} \left| \left\langle f, \sum_{i=1}^{N}{a_i \delta_{x_i}} \right\rangle  \right|^2 \\
&= \sup_{f \in V \atop \|f\|_{L^2} \leq 1} \left| \left\langle f, \pi \sum_{i=1}^{N}{a_i \delta_{x_i}} \right\rangle  \right|^2\\
&= \left\| \pi \sum_{i=1}^{N}{a_i \delta_{x_i}} \right\|^{2} = \mathcal{E}_{\leq s}.
\end{aligned}
\end{equation}
This naturally motivates $\mathcal{E}_{\leq s}$ as an interesting objective function that one should try to minimize. Different function classes defined via weighted Fourier coefficients, for example Sobolev spaces, would yield slightly different terms that might motivate different functionals; we consider this an interesting direction for further research.
\subsubsection{Results for a single annealing run}
\begin{center}
\begin{figure}[ht]
 \includegraphics[width=5.0in]{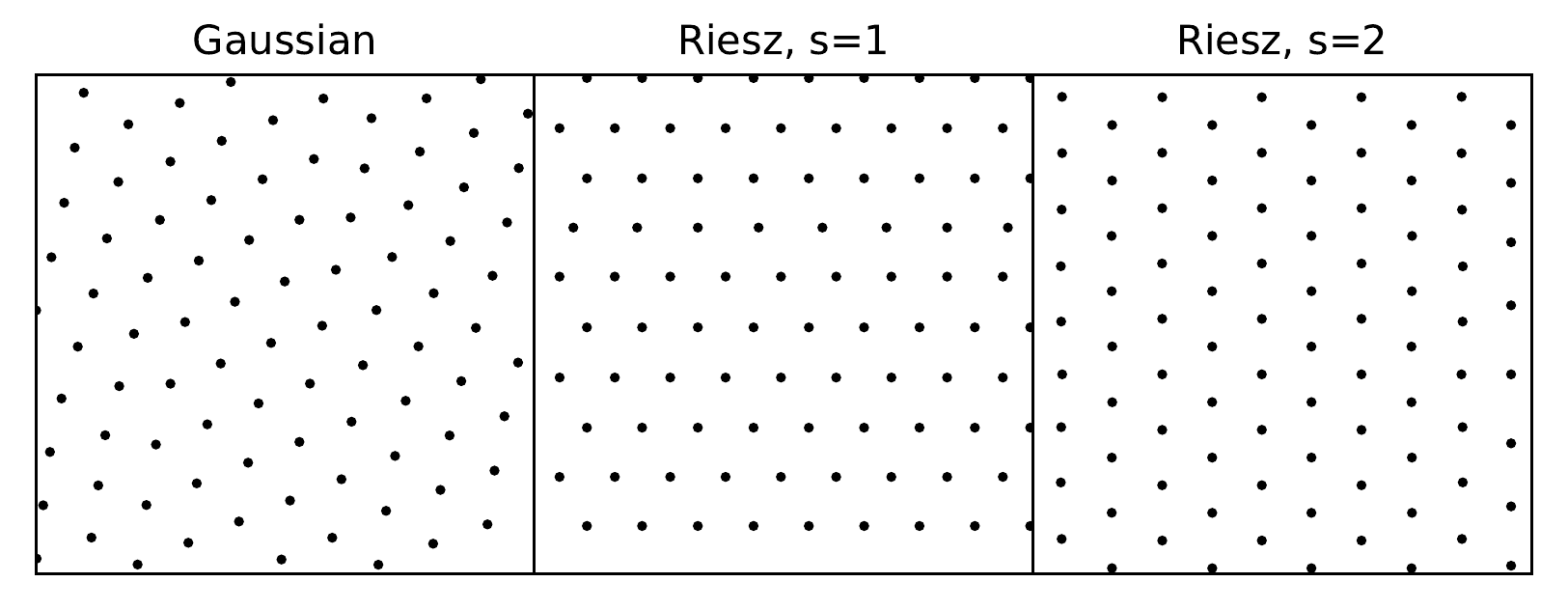}
\caption{Minimum energy configurations of $N=89$ support points in the case of the standard 2-torus for the energy functionals $E_{\rm Gaussian}$, and $E_{{\rm Riesz}, s},s=1,2$, respectively. }\label{fig:config:demo}
\end{figure}
\end{center}

\begin{center}
\begin{figure}[ht]
 \includegraphics[width=5.0in]{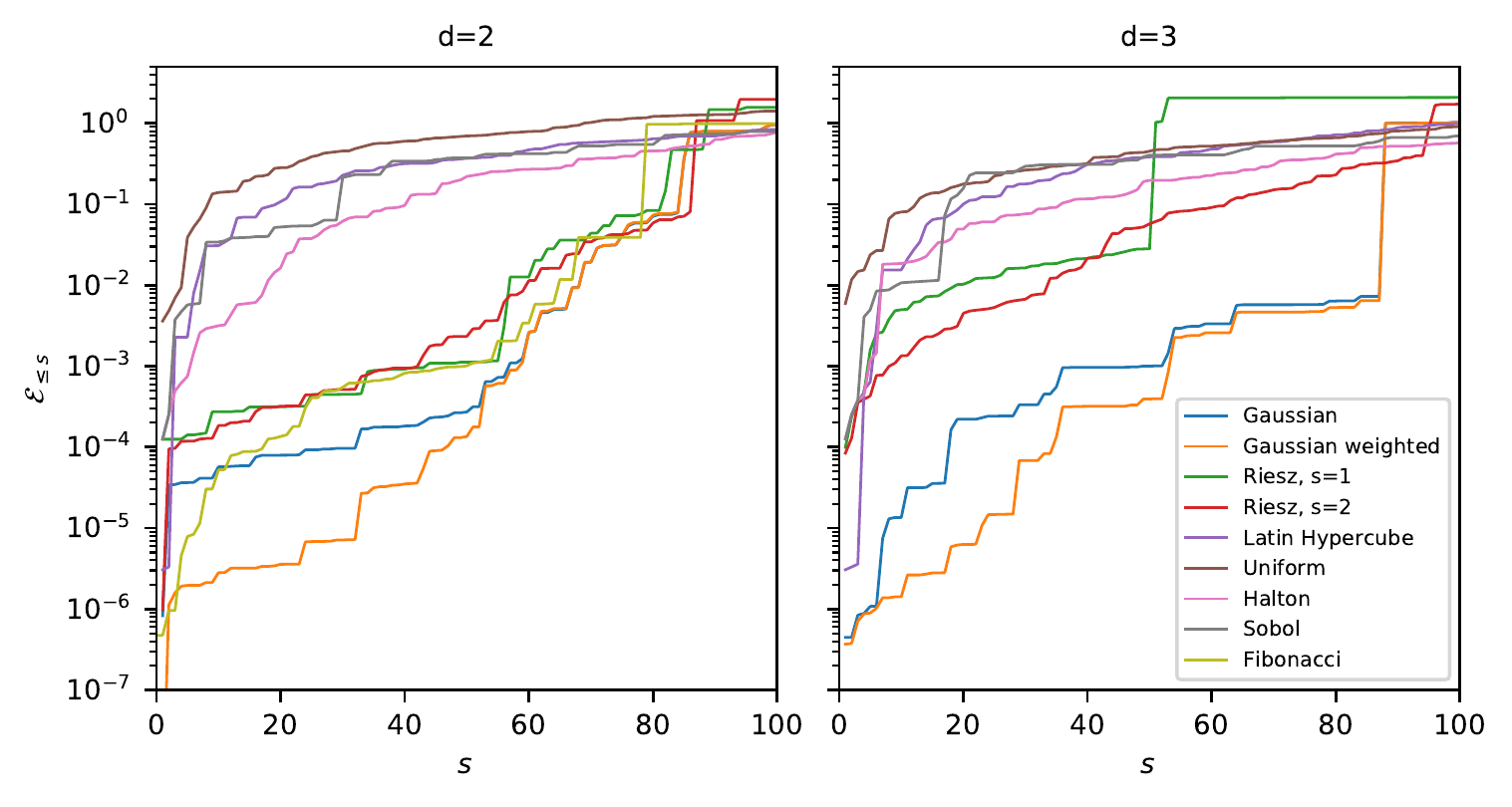}
\caption{Quadrature errors for $N=89$ points on $\mathbb{T}^{d},~d=2,3$. (color figure online)}\label{fig:torus:89}
\end{figure}
\end{center}

\begin{center}
\begin{figure}[ht]
 \includegraphics[width=5.0in]{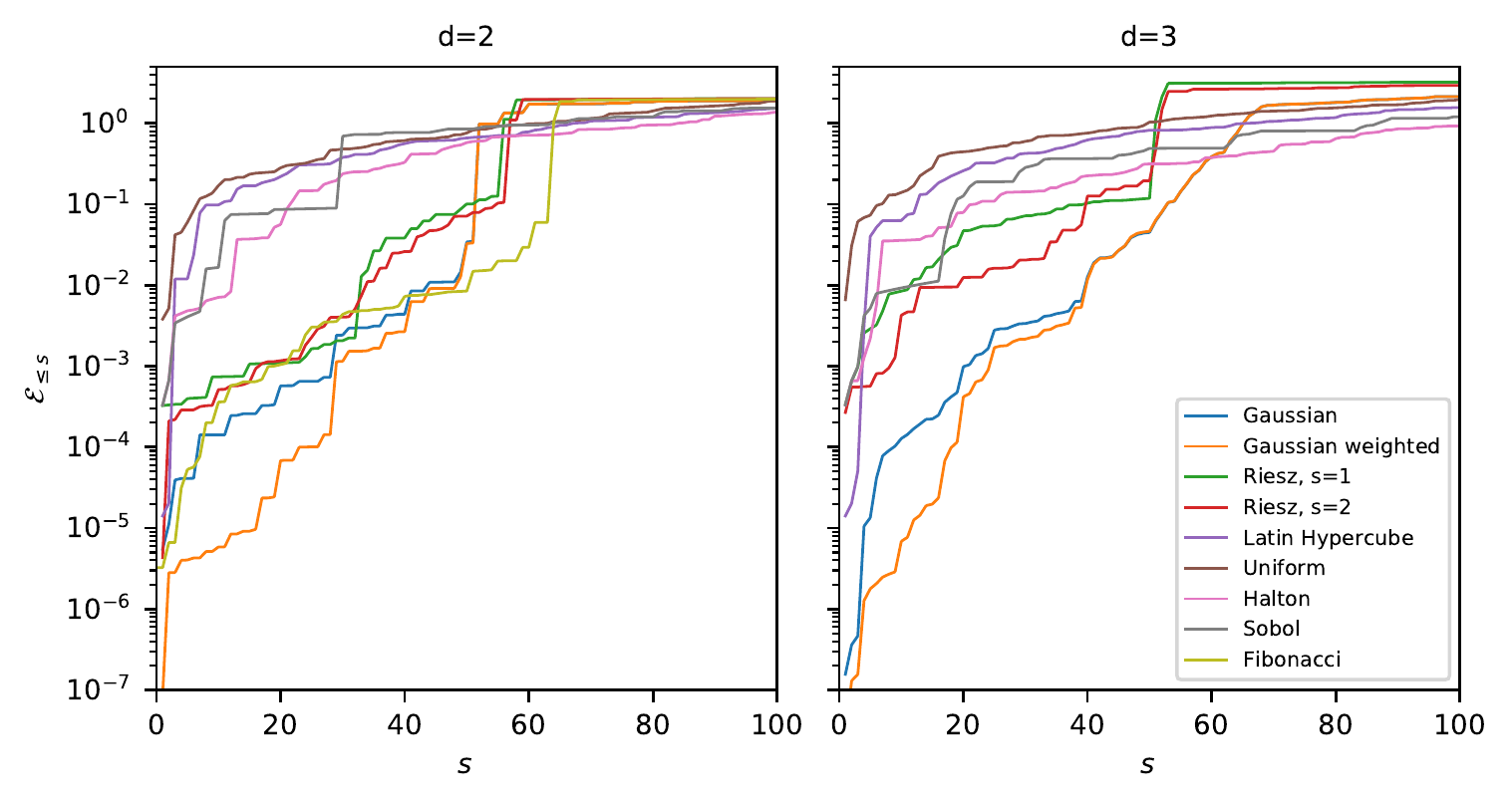}
\caption{Quadrature errors for $N=55$ points on $\mathbb{T}^{d},~d=2,3$. (color figure online)}\label{fig:torus:55}
\end{figure}
\end{center}

We consider setups of $N=89$ and $N=55$ quadrature points in the 2-torus and on the 3-torus. \cref{fig:config:demo} shows examples of point sets of  $N=89$ support points generated by minimizing the Gaussian energy functional $E_{\rm Gaussian}$, and the Riesz energy functionals $ E_{{\rm Riesz}, s}, s=1,2$, respectively. Choosing the number of quadrature points to be Fibonacci numbers allows us to compare against Fibonacci lattices on $\mathbb{T}^{2}$. Fibonacci point sets are constructed as an integration lattice (see e.g. \cite{owen1}) and satisfy the criteria of a good lattice point set; see \cite{niederreiter}. They have been shown to be optimal in certain error measures; see e.g. \cite{zaremba}.
Figure \ref{fig:torus:89} and Figure \ref{fig:torus:55} show the quadrature error for $N=89$ and $N=55$ quadrature points, respectively. In terms of the considered error measure we find that point sets constructed by our approach, i.e., as minimal energy configurations of $E_{\rm Gaussian}$,
compare very favourably to existing methods. The incorporation of weights leads to a further reduction of the integration error to up to more than one order of magnitude in the low frequency range (i.e. small values of the order index $s$) and only the Fibonacci lattices result in a comparable integration error. Other classical QMC point sets, namely the Sobol point sets and the Halton point sets as well as the considered (stratified) Monte-Carlo point sets, are clearly outperformed in the low frequency range, i.e.,  for small values of the order index $s$ we find in all considered setups that the error measure $\mathcal{E}_{\leq s}$ is several magnitudes smaller for the weighted points constructed  by our approach  in comparison to the above mentioned point sets. Similarly, we find for the considered Riesz energy functionals that the corresponding point sets result in an integration error which is significantly worse than the point sets constructed as minimizers of the Gaussian energy functional. Moreover, \cref{fig:torus:89} and Figure \ref{fig:torus:55} suggest that there is a tradeoff between reducing the integration error for small values of the order index $s$ versus reducing the integration error for high values of the order index $s$, i.e.,  while as described above the minimum energy point sets perform particularly well for small order indices, we find that for large $s$ (e.g. $s\geq 60$ in the case of $N=55, d=2$) the integration error of the minimum energy point sets is consistently worse across all considered setups than the integration error of the considered QMC point sets and MC point sets.

\subsubsection{Results for multiple annealing runs}
\begin{center}
\begin{figure}[ht]
 \includegraphics[width=5.0in]{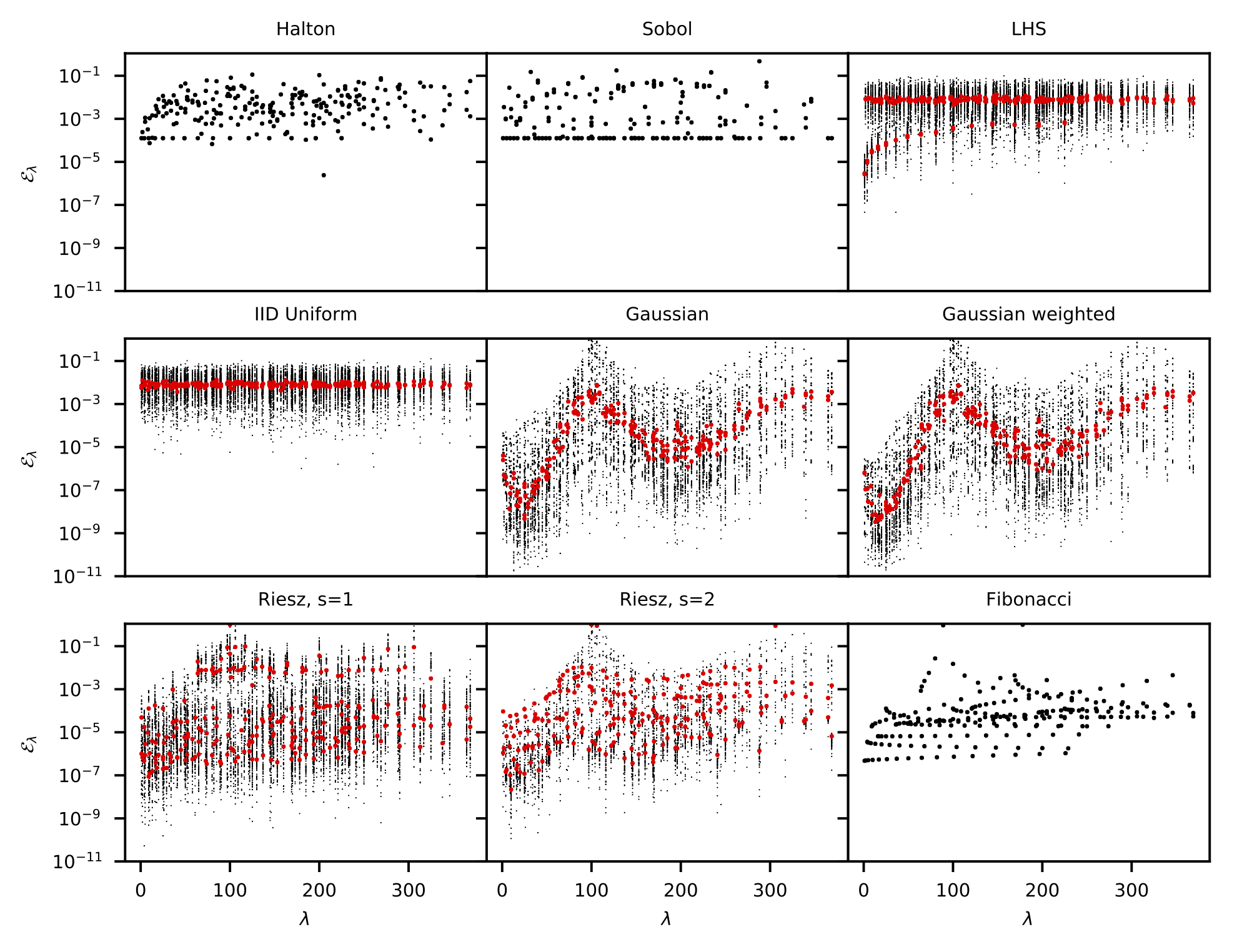}
\caption{Scatter plot of $\lambda \in \sigma(-\Delta)$ vs. integration error $\mathcal{E}_{\lambda}$ for $N=89$ quadrature points on the $\mathbb{T}^2$. For non-deterministic methods the resulting error values of $N_{r}=50$ independent samples are superimposed with the median values marked in red. (color figure online)}\label{fig:torus2:boxplot}
\end{figure}
\end{center}

\begin{center}
\begin{figure}[ht]
 \includegraphics[width=5.0in]{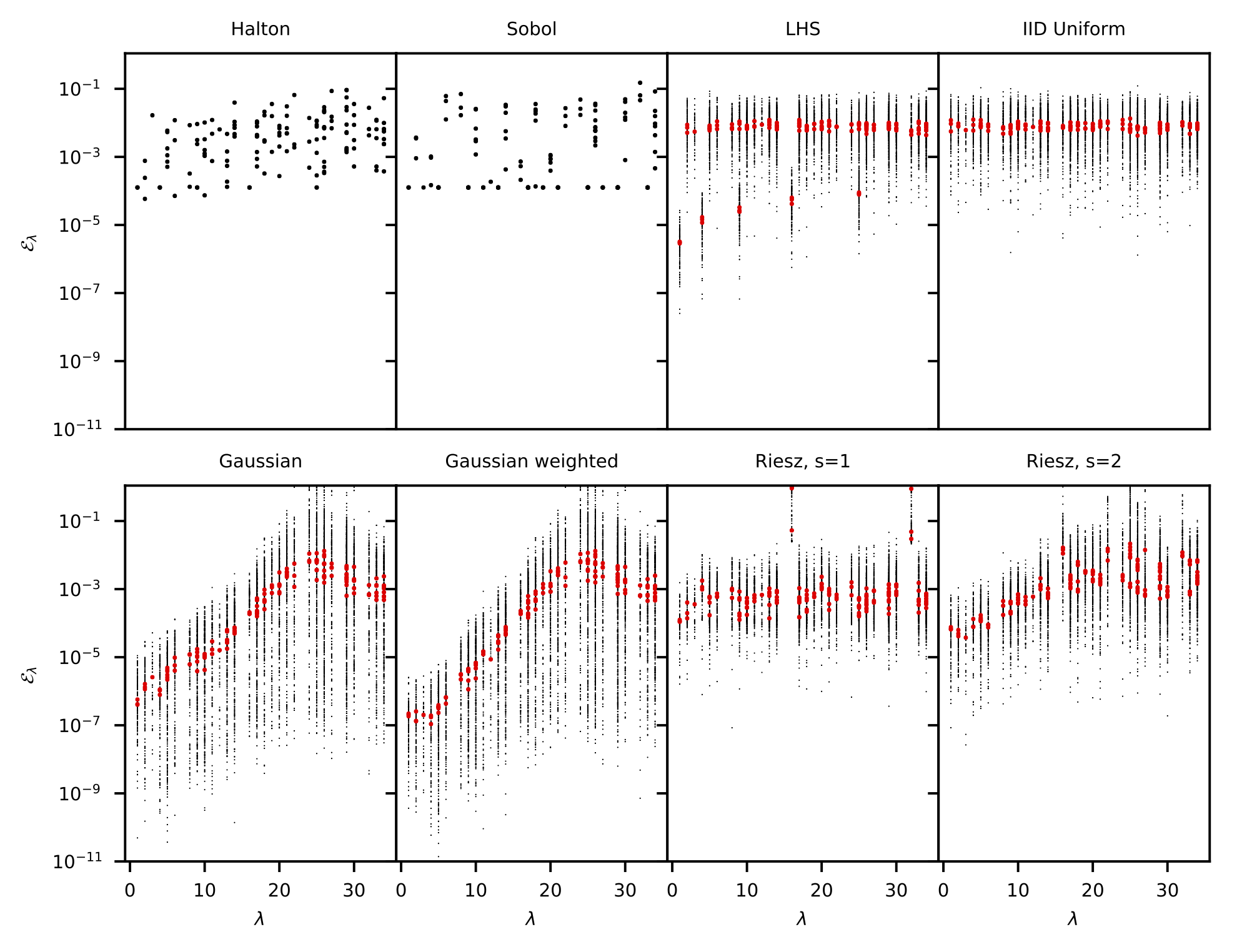}
\caption{Scatter plot of $\lambda \in \sigma(-\Delta)$ vs. integration error $\mathcal{E}_{\lambda}$ for $N=89$ quadrature points on the 3-torus. In the case of non-deterministic methods we consider $N_{r}=50$ independent samples. The format of the figure is the same as the format of \cref{fig:torus2:boxplot}. (color figure online)}\label{fig:torus3:boxplot}
\end{figure}
\end{center}

Since the Riesz minimum energy point sets and the Gaussian miminum energy point sets are constructed using a non-deterministic procedure one can expect fluctuations in the quality of the resulting quadrature points. More generally, since the corresponding global minimisation problem is intractable for a sufficiently high number of support points, one will in practice always end up with point sets which correspond to non global local minima of the respective energy functional. For this reason it is sensible to not only compare single point sets, but instead assess the integration error for an ensemble of point sets obtained from independent annealing runs. \cref{fig:torus2:boxplot} and \cref{fig:torus3:boxplot} show the resulting integration errors measured in terms of
$
\mathcal{E}_{\lambda}, \lambda \in \sigma(-\Delta),
$
for $N=89$ quadrature points for point sets obtained from  $N_{r}=50$ independent annealing runs applied to the energy functionals considered in the previous section. For both the 2-torus (see \cref{fig:torus2:boxplot}) and the 3-torus (see  \cref{fig:torus3:boxplot}), we find that the qualitative behaviour of the integration error measured in terms of $\mathcal{E}_{\lambda}$ is consistent with what we observed in the preceding section for the considered methods. In particular, the integration error $\mathcal{E}_{\lambda}$ of Gaussian minimum energy points appears to behave cyclic as $\lambda$ increases with comparably small integration errors for small values of $\lambda$ and comparably large integration errors for values $\lambda \approx \widetilde{\lambda}$ with $\widetilde{\lambda}=106$ in the case of $M=\mathbb{T}^{2}$ and $ \widetilde{\lambda} = 26$ for $M=\mathbb{T}^{3}$.
Although less pronounced, a similar qualitative behaviour can be also found for the integration error of the Riesz minimum energy point sets and the Fibonacci point set, whereas in the case of the other QMC point sets and the MC  point sets there is no obvious relationship between the magnitude of $\lambda$ and the integration error for the corresponding eigenfunction. This explains the crossing as $s\rightarrow \infty$ of the curves of $\mathcal{E}_{\leq s}$ associated with different points sets in \cref{fig:torus:89}. More specifically, for Gaussian minimum energy points (both weighted and unweighted), we observe that the median value of the integration errors for eigenfunctions associated with eigenvalues $\lambda <\widehat{\lambda}$ (with $\widehat{\lambda}\approx80$ for $M=\mathbb{T}^{2}$ and  $\widehat{\lambda}\approx15$ for $M=\mathbb{T}^{3}$) is up to three magnitudes smaller than the integration error of the considered MC and most QMC methods. Even when considering the maximum of $\mathcal{E}_{\lambda}$ across all 50 annealing runs, the Gaussian minimum energy point sets compare favourably for $\lambda \leq \widehat{\lambda}$ to other methods. Only the Fibonacci point set on $M=\mathbb{T}^{2}$ results in integration errors which are of comparable magnitude for $\lambda \leq \widehat{\lambda}$. 





\subsection{On the sphere.}
We shall compare our approach against established Quasi-Monte Carlo point sets for integration on the sphere, i.e.,
\begin{enumerate}
\item Spherical Designs \cite{delsarte}
\item Fibonacci numerical integration on a sphere (see \cite{hannay})
\item Minimal Energy Configurations for the Riesz energies on the sphere
\end{enumerate}
as well as the i.i.d. uniform random point set on the sphere. The spherical design quadrature points against which we compare were obtained from \cite{womersleyurl}.
\subsubsection{Construction of quadrature points and weights}
We construct minimum energy configuration for the considered energy functionals using a simulated annealing scheme based on a Langevin diffusion constrained on $\mathbb{S}^2$. Weights for point sets corresponding to local minima of the Gaussian energy functional are constructed in the same way as described in \cref{sec:torus:construction}, and we note that as for the point sets considered for the 2-torus and 3-torus, the closed form solution \eqref{eq:sol:weights} of the associated optimisation is positive. 
In order to simplify computations we approximate the geodesic distance metric on $\mathbb{S}^{2}$  by the Euclidian distance in $\mathbb{R}^{3}$  in the computation of the energy functionals $\widetilde{E}_{\rm Gaussian}$, and $E_{{\rm Riesz}, s}$. As explained in Section~\ref{sec:torus:construction}, at least for the Gaussian energy functional  this approximation is well justified since tails of the Gaussian kernel decay sufficiently fast so that the energy contribution of pairs of points for which this local approximation of the geodesic distance does not hold, are negligible.  
\subsubsection{Evaluation of quadrature error}
As for the torus case, we evaluate the error in terms of the integration error incurred for eigenfunctions of the Laplace operator, which on $\mathbb{S}^{2}$ are given by the spherical polynomials. Let $\left ( a_{j},x_{j} \right )_{1 \leq j \leq N}$ be prescribed (weighted) point set. In analogy to the error measure considered above for the $d$-torus, we define
\begin{equation*}
\mathcal{E}_{\lambda_{s}} = \abs*{\frac{1}{N} \sum_{j=1}^{N}  a_{j}Y^{m}_{l}(x_{j})}^{2}
\end{equation*}
where $Y^{m}_{l} : \mathbb{S}^{2} \rightarrow \mathbb{C}$ denote the spherical polynomial of degree $l\in \mathbb{N}$ and order $m$, i.e.,
\[
-\Delta Y^{m}_{l} = l (l+1) Y^{m}_{l}= \lambda_{s}Y^{m}_{l}, ~~ -l \leq m \leq l.
\]
As in the torus case we denote the cumulative sum of $\mathcal{E}_{\lambda_{s}}, s=1,2,\dots$, as $\mathcal{E}_{\leq s}$, i.e.,
\begin{equation*}
\mathcal{E}_{\leq s} := \sum_{l=1}^{s}\mathcal{E}_{\lambda_{l}}.
\end{equation*}
\subsubsection{Numerical results}
\begin{center}
\begin{figure}[ht]
 \includegraphics[width=5.0in]{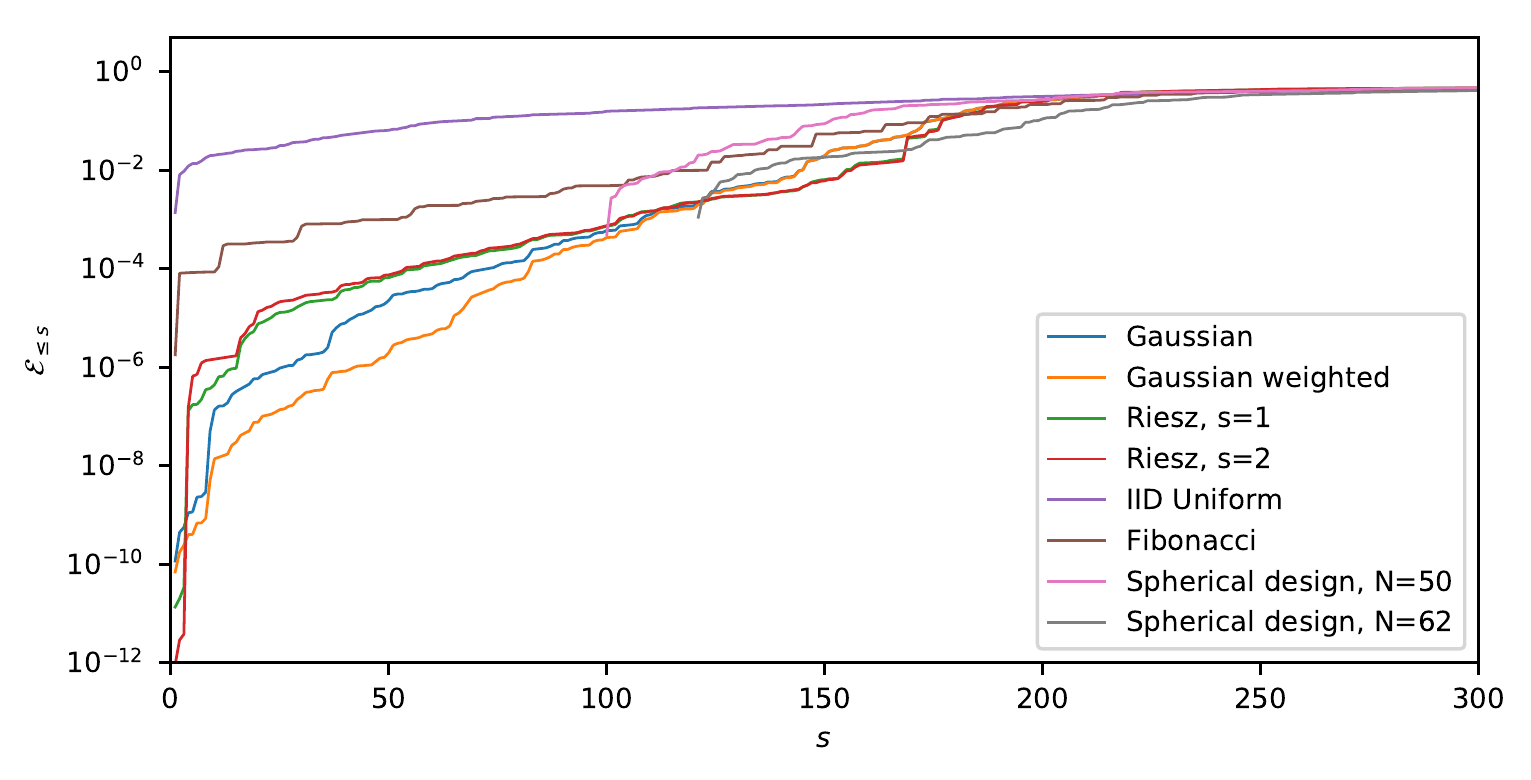}
\caption{Quadrature errors for $N=55$ points on $\mathbb{S}^{2}$. (color figure online)}\label{fig:sphere:2}
\end{figure}
\end{center}

\begin{center}
\begin{figure}[ht]
 \includegraphics[width=5.0in]{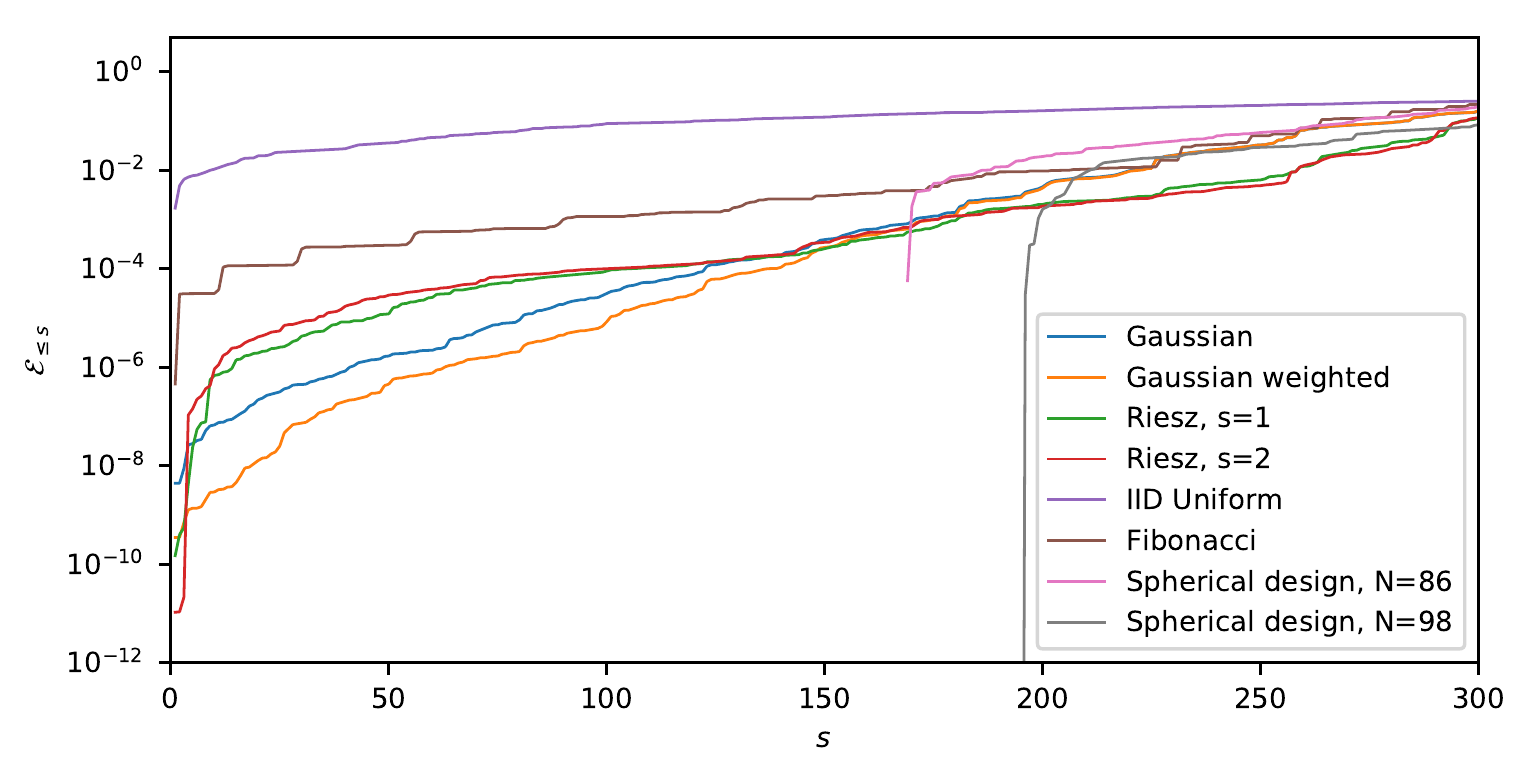}
\caption{Quadrature errors for $N=89$ points on $\mathbb{S}^{2}$. (color figure online)}\label{fig:sphere:1}
\end{figure}
\end{center}

We consider the case of $N=89$ and $N=55$ quadrature points on the sphere.  \cref{fig:spheres} (A) shows a point set of $N=89$ quadrature points generated by our approach. Since to the best of our knowledge there are no spherical design point sets known for these numbers of points, we compare against spherical design point sets of size $86$ and $96$ in the case of $N=89$ and against spherical design point sets comprising $50$ and $62$ points in the case of $N=55$. By construction these point sets allow exact integration of spherical harmonics up to degrees $l= 9,10,12,13$ for the point set sizes $50,62,86,99$, respectively. \cref{fig:sphere:1} and \cref{fig:sphere:2} show the cumulative integration $\mathcal{E}_{\leq s}$ error for $N=89$ and $N=55$, respectively.


\begin{center}
\begin{figure}[ht]
 \includegraphics[width=5.0in]{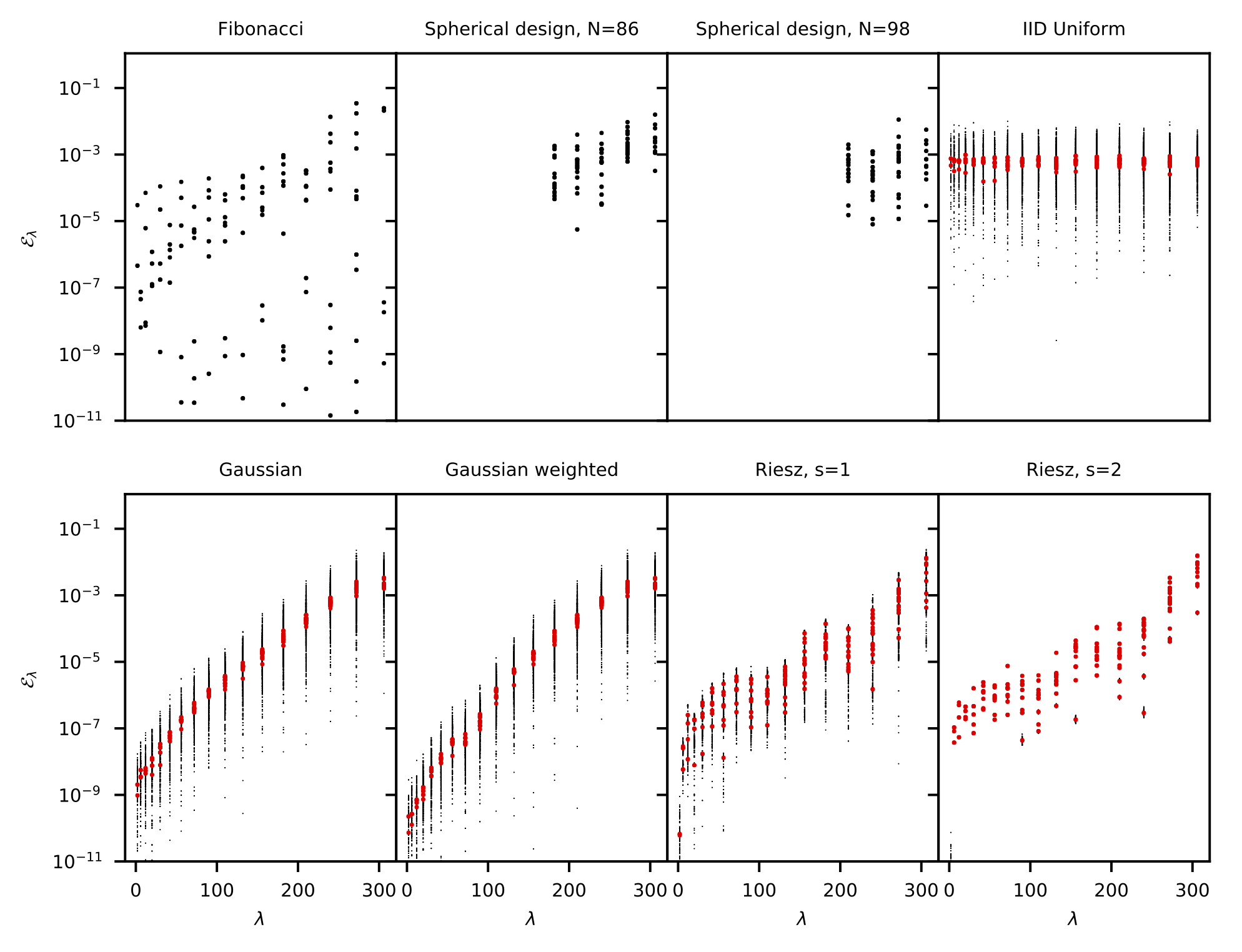}
\caption{Scatter plot of $\lambda \in \sigma(-\Delta)$ vs. integration error $\mathcal{E}_{\lambda}$ for $N=89$ quadrature points on the 2-sphere. For non-deterministic methods the resulting error values of $N_{r}=50$ independent sample are superimposed with median values marked in red. (color figure online)}\label{fig:sphere:boxplot}
\end{figure}
\end{center}

 In both setups the Fibonacci point sets and the i.i.d. uniform point set are clearly outperformed by the Gaussian minimum energy point sets in terms of the considered error measure. Interestingly,  for Riesz minimum energy point sets outperform the Gaussian minimum energy point sets in terms of the error measure $\mathcal{E}_{\leq s}$ in the cases of $s\leq 2$ or $s \geq 150$ for $N=89$ support points, and $s\geq 2$ or $s \geq 100$ for $N=55$ support points.
This observation is supported by \cref{fig:sphere:boxplot} where the integration errors $\mathcal{E}_{\lambda}, \lambda \in \sigma(-\Delta)$ of points sets obtained from and $N_{r} = 50$ independent annealing runs are superimposed in a scatter plot of the same format as the above \cref{fig:torus2:boxplot} and \cref{fig:torus3:boxplot}.
In all the considered setups our approach is outperformed by the spherical design point sets. We emphasize, however, that these spherical design point sets are built so as to minimize the error functional, are generally very hard to construct, and may only exist for certain values of $N$.
For general smooth manifolds the explicit construction of point sets which exactly integrate a prescribed number of  Laplacian eigenfunctions is a not well understood problem (analogous notions of designs seem to exist in a fairly abstract setting, see \cite{steingraph} and \cite{sey}, but we are not aware of any reliable way they could be constructed on manifolds that are not spheres).
\FloatBarrier
\subsection{On other manifolds.}
\begin{figure}
    \centering
    \begin{subfigure}[b]{0.5\textwidth}
        \includegraphics[width=\textwidth]{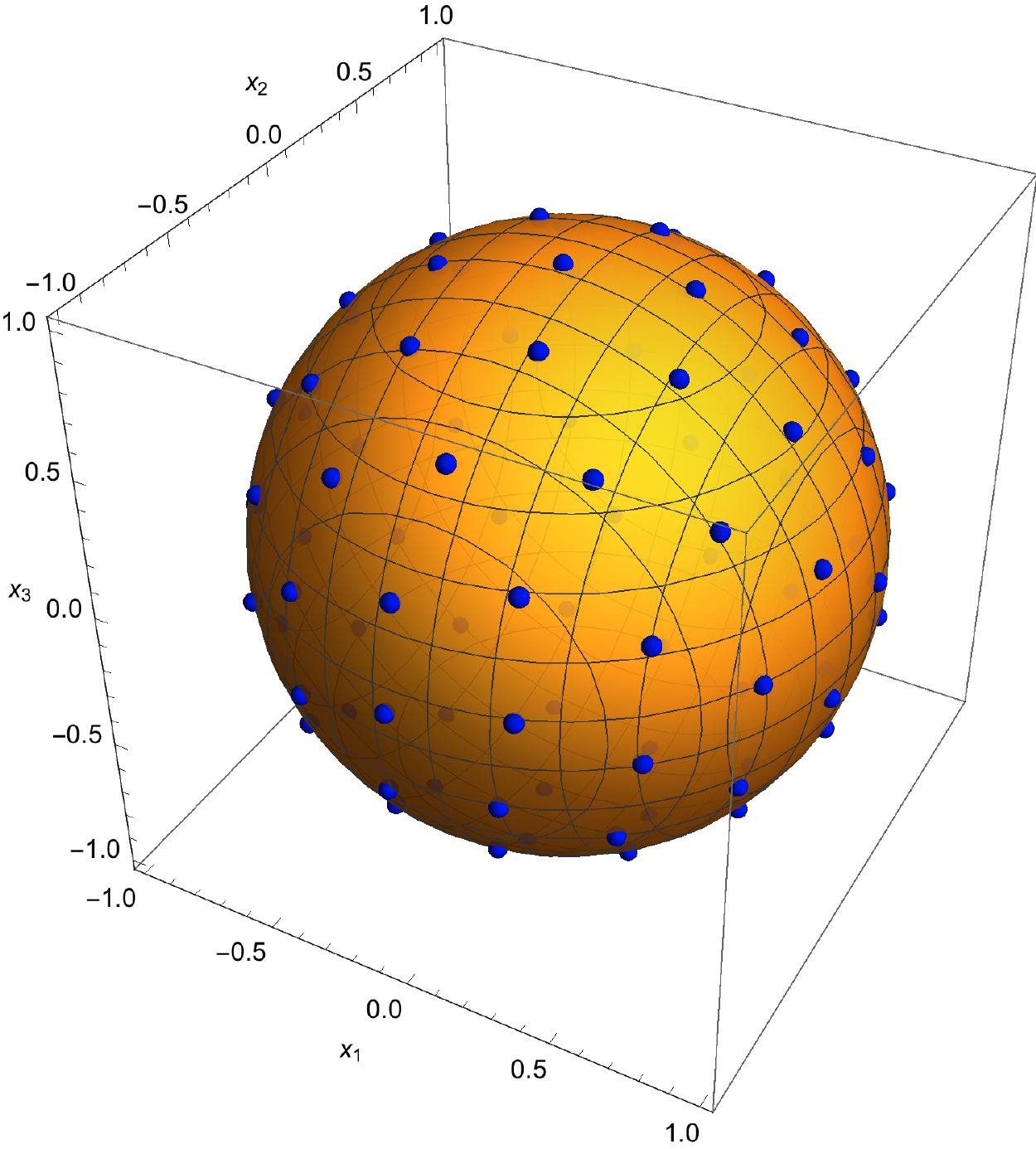}
        \caption{}
        \label{fig:sphere}
    \end{subfigure}
    ~ 
    \begin{subfigure}[b]{0.5\textwidth}
        \includegraphics[width=\textwidth]{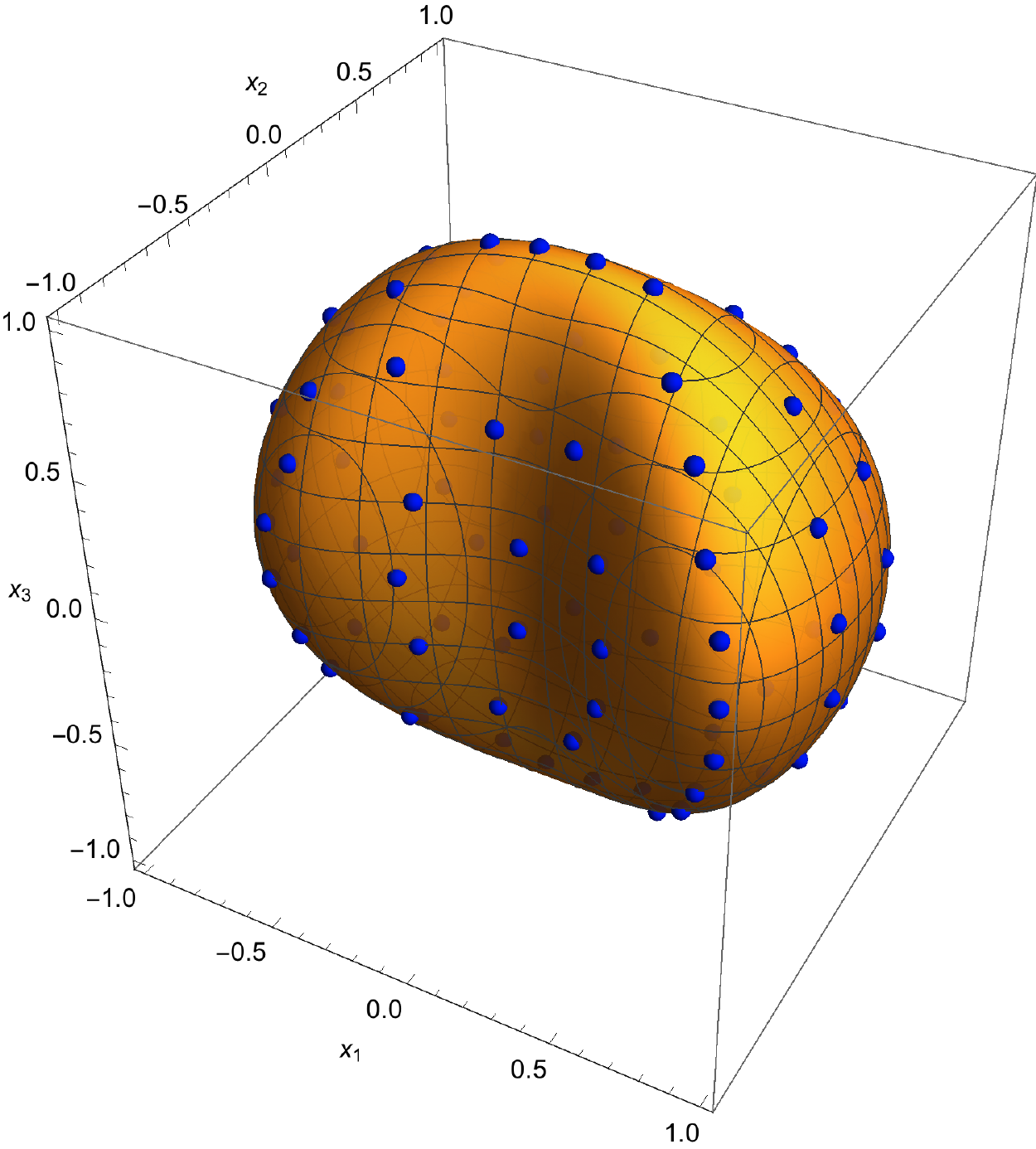}
        \caption{}
        \label{fig:dented:sphere}
    \end{subfigure}
    ~ 
    \caption{Support points constructed as minimum energy configurations of the energy functional $E_{\rm Gaussian}$ on the sphere (A), and on a ``dented'' sphere (B) corresponding to the level set \eqref{eq:dented:sphere} with $\alpha=1/10$. (color figure online)}\label{fig:spheres}
\end{figure}

In order to demonstrate the applicability of the proposed approach to manifolds beyond the ones described above, we consider two manifolds for which the construction of support points is a less well studied problem. Unlike in the cases of $M= \mathbb{T}^{d}$, and $M=\mathbb{S}^{d-1}$, there are no closed form solutions for the Laplacian eigenfunctions for the manifolds considered in this section. We thus provide only qualitative results. We first consider a manifold which embedded in the 3-dimensional Euclidian space corresponds to the solution set of
\begin{equation}\label{eq:dented:sphere}
x_1^2+\frac{x_2^2}{\alpha+x_1^2}+x_3^2 = 1,
\end{equation}
where $\alpha$ is a positive scalar. This manifold is an element of the same homology group as $\mathbb{S}^{2}$, and we refer to it as a ``dented'' sphere. The magnitude of the scalar $\alpha$ determines the strength of the dent. \cref{fig:spheres} (B) shows a point set constructed with our approach in the case of $\alpha=1/10$ and $N=89$ support points. 
\\

As a second example we consider the construction of support points on a compactified version of the Poincar\'e disk model (see e.g. \cite{anderson2006hyperbolic}). More specifically, we consider the Poincar\'e disk model constrained to the subset which in Euclidian space corresponds to the disk centred in the origin with radius $r=4/5$. Strictly speaking this model does not fall into the class of manifolds specified in the exposition of our approach in the introduction of this paper. However, assuming Neumann boundary conditions the spectrum of the Laplace-Bletrami operator can be shown to be discrete, and the associated heat propagator is self-adjoint, thus the identity \eqref{eq:identity:1} and the inequality \eqref{eq:inequality:1} remain valid. This justifies the application of our approach.  
The Poincar\'e disk model is related to the hyperboloid model, i.e., the manifold associated with the solution set of the equation
\begin{equation}\label{eq:hyperboloid}
x_{1}^{2}+x_{2}^{2}-x_{3}^{2} = -1,
\end{equation}
by the fact that the projection \footnote{Geometrically, this projection can be interpreted as the intersection of the line connecting the point $(x_{1},x_{2},x_{3})$ and the point $(0,0,-1)$ with the hyperplane spanned by the first two canonical basis vectors in $\mathbb{R}^{3}$.}
\begin{equation}\label{eq:projection}
(x_{1},x_{2},x_{3}) \mapsto   \left ( \frac{x_{1}}{1 + x_{3}} , \frac{x_{2}}{1 + x_{3}} \right ) 
\end{equation}
maps geodesics of the upper sheet of the hyperboloid model onto geodesics of the Poincar\'e disk model. 
We use this relation to construct support points on the compactified Poincar\'e disk model by first minimizing the energy functional $E_{\rm Gaussian}$ on the upper sheet of a suitably compactified version of the hyperboloid model which is the manifold associated with the solution set of \eqref{eq:hyperboloid} with the additional constraint,
\begin{equation}\label{eq:add:constraint}
x_{3} \leq  \frac{ 1+ r^{2}}{1 - r^{2}},
\end{equation}
where as above $r=4/5$. The construction of support points in this way is illustrated in \cref{fig:poincare} (A), which  shows a point set of $N=150$ points on the compactified upper sheet of the hyperboloid constructed as a minimum energy configuration of the energy functional $E_{\rm Gaussian}$. The image of this point set under the projection \eqref{eq:projection}, is shown in \cref{fig:poincare} (B). As  explained above this point set corresponds by construction to a minimum energy configuration of the energy functional $E_{\rm Gaussian}$ on the compactified  Poincar\'e disk model.




\begin{figure}
    \centering
    \begin{subfigure}[b]{0.5\textwidth}
        \includegraphics[width=\textwidth]{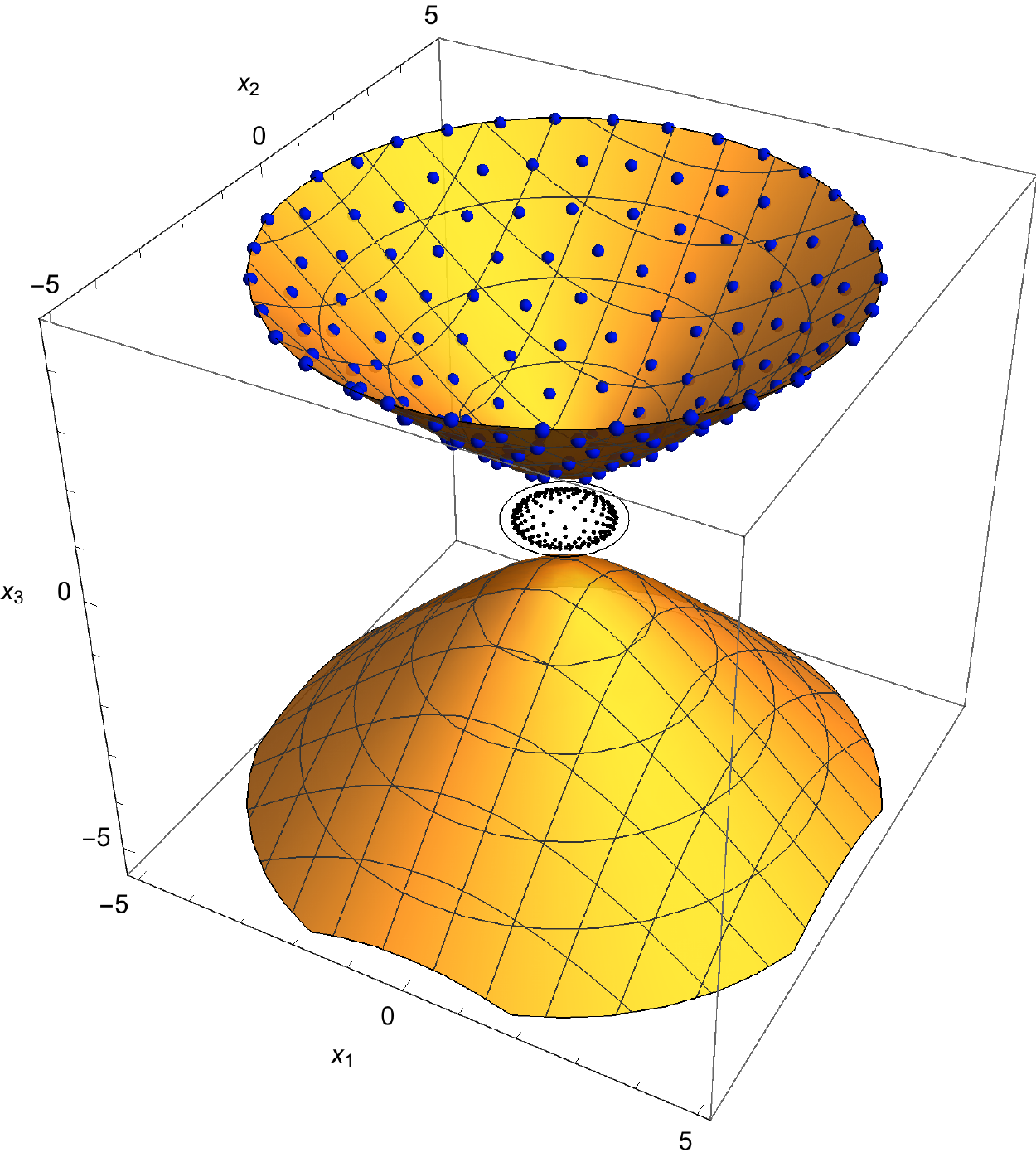}
        \caption{}
        \label{fig:sphere}
    \end{subfigure}
    ~ 
    \begin{subfigure}[b]{0.5\textwidth}
        \includegraphics[width=\textwidth]{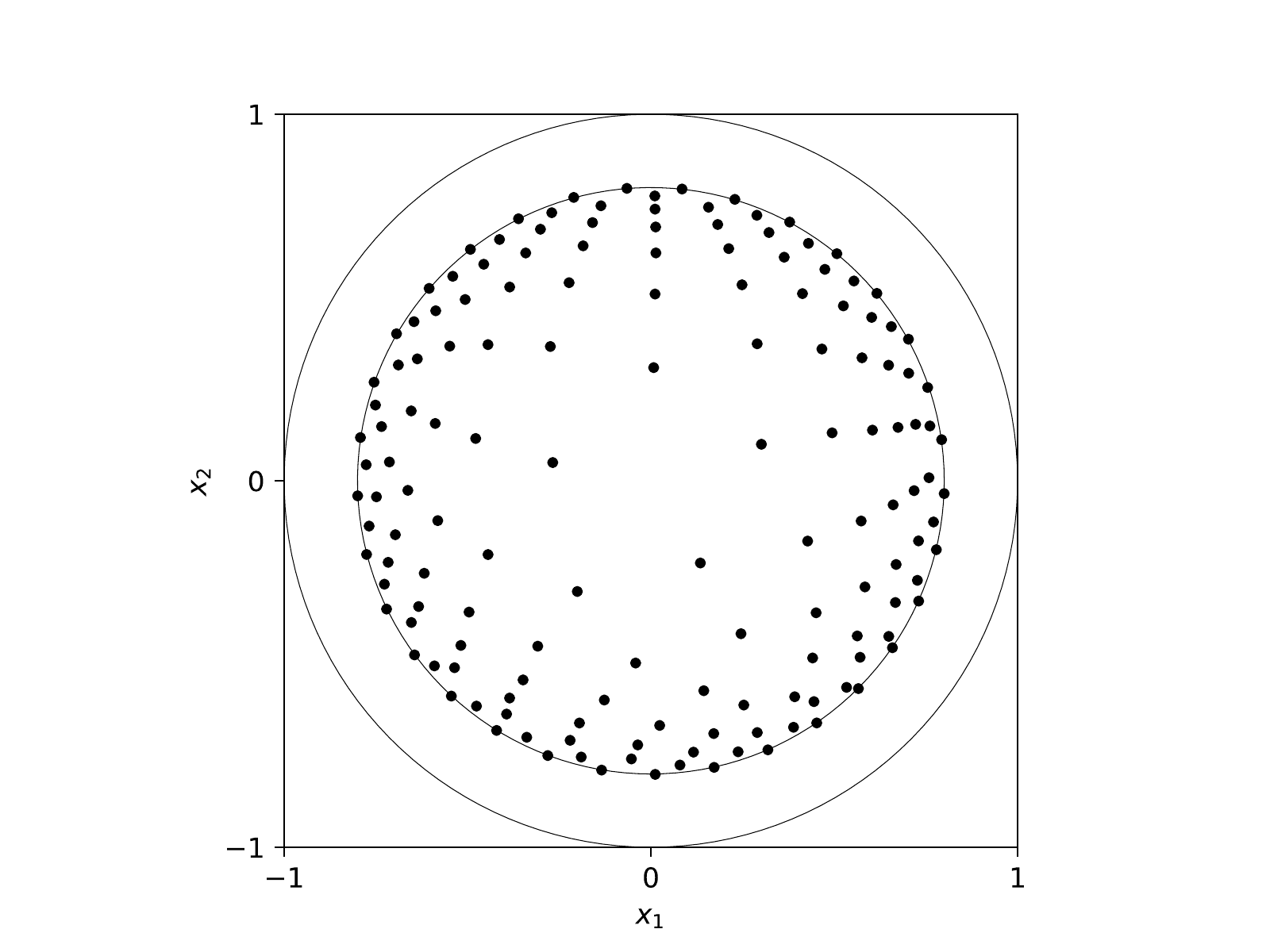}
        \caption{}
        \label{fig:poincare}
    \end{subfigure}
    ~ 
    \caption{Support points constructed as minimum energy configuration of the energy functional $E_{\rm Gaussian}$ on a compactified version of the upper hyperboloid (A), and the projection of these onto the unit disk (B) resulting in support points for the compactified Poincar\'e disk model. (color figure online)}\label{fig:poincare}
\end{figure}

\FloatBarrier
\section{Conclusion}
It is naturally desirable to have quadrature rules that perform well on low-frequency Laplacian eigenfunctions since those are the smoothest orthogonal functions on any geometry. We have shown that using this as a starting point, we obtain a natural functional on the set of weighted points 
$$ \sum_{i,j =1}^{N}{ a_i a_j p(t,x_i, x_j) } \rightarrow \min  \qquad \mbox{where}~t \sim N^{-\frac{2}{d}}.$$
Minimizers of this functionals have, by construction, a very small weighted error over the first few Laplacian eigenfunctions. This is difficult to use in practice since the heat kernel
is typically not available. However, since $t \sim N^{-2/d} \ll 1$, one would expect very similar results when replacing the heat kernel with Varadhan's short time asymptotic
$$ \sum_{i,j =1}^{N}{ a_i a_j \exp\left(-\frac{d(x_i,x_j)^2}{4t}\right) } \rightarrow \min  \qquad \mbox{where}~t \sim N^{-\frac{2}{d}}.$$
This produces a geometry-independent functional on weighted points that, by design, is expected to have a small error on low-frequency functions. A nice by-product is that
it also automatically generates weights. We have shown the method to be competitive against methods that have been explicitly tailored for the sphere and the torus. We have not compared the method to methods on other manifolds since there are not that many universal design rules (though this is becoming a more popular topic of research, see \cite{steinpoints, supportpoints}).

Our method is superior to points minimizing the Riesz energy; this is, considering the singularities of the  Fourier transform of the Riesz kernel, perhaps not surprising.

There are many open questions that remain.

\begin{enumerate}
\item How does the method do on general manifolds? What would be reasonable point sets to compare against? Is there a version on graphs $G=(V,E)$?
\item Are there ways to speed up the computation? Is the quality of the point sets strongly depending on the numerical algorithm with which it was obtained?
\item It was recently pointed out by Ehler, Graef \& Oates \cite{ehler} that adding weights can substantially improve the integration error in reproducing kernel Hilbert spaces. Our method naturally provides a way to obtain weights for any given set of points, it would be interesting to understand to which extent these weights have good properties.
\item Are these points especially useful for the discretization of parabolic partial differential equations since the bulk of their dynamical behavior is at low frequencies? (This question is due to Manas Rachh.)
\item Instead of the heat kernel, one could take other kernels with slightly slower decay in the spectrum; this would possibly lead to slightly larger errors on
low-frequency eigenfunctions but could possibly yield some improvement at the eigenfunctions after the first $\sim N$ eigenfunctions. The Riesz kernel on
the torus is such an example but we are not aware of the corresponding spectral theory for the Riesz kernel on general manifolds.
\end{enumerate}

\FloatBarrier

\appendix
\section{Simulated annealing scheme for $M=\mathbb{T}^{d}$}\label{ap:sec:annealing:torus}
In oder to find good approximations of global minimizers of the energy functions  $ E_{\rm Gaussian}$ and $E_{{\rm Riesz},s}$,  we use an annealing scheme based on the stochastic dynamics described by an It\^o-diffusion of the form 
\begin{equation}\label{eq:ULD}
\begin{aligned}
\dd \x(t) & = \p(t) \dd t,\\
\dd \p(t) &= - \nabla_{\x}U(\x(t)) \dd t - \gamma \p(t) \dd t + \sqrt{2 \gamma\beta^{-1}(t)} \dd\W(t),\\
\end{aligned}
\end{equation}
where 
\begin{enumerate}
\item $\W = (W_{1}, \dots, W_{Nd})$, with $W_{i}, 1\leq i \leq Nd$ being independent Wiener processes,
\item $\x = (\x_{1},\dots,\x_{N})$ so that $\x_{i}(t) \in \mathbb{T}^{d},1\leq i \leq N$,
\item $\gamma>0$, and $\beta \in \mathcal{C}([0,\infty), \mathbb{R}_{+} )$,
\end{enumerate}
with $U \in \left \{ E_{\rm Gaussian}, E_{{\rm Riesz},s} \right \}$.
The stochastic differential equation (SDE) \eqref{eq:ULD} is known in the statistical physics literature as the {\em underdamped Langevin equation}. The underdamped Langevin equation can be viewed as a stochastic perturbed version of Hamilton's equation associated with the Hamiltonian $ H(\x,\p)= U(\x) +\frac{1}{2}\norm{\p}_{2}^{2}$. The remaining terms in \eqref{eq:ULD} model the exchange of energy with a heat bath; see e.g. \cite{pavliotis} for more details.  The parameter $\gamma>0$ determines the strength of the coupling and as such can be interpreted as a friction coefficient. $\beta(t)>0$ can be interpreted as the inverse temperature of the heat bath, and the function $\beta : [0,\infty) \rightarrow [0,\infty)$ is commonly referred to as a {\em cooling schedule}.\\
We discretize  \eqref{eq:ULD}  using the well studied ``BAOAB''-splitting scheme (see \cref{alg:annealing}), which as a symmetric stochastic splitting scheme is of weak second order accuracy in the discretization/stepsize parameter $\deltat$; see \cite{leimkuhler1,leimkuhler2} for details. We parametrize \cref{alg:annealing}  with a cooling schedule $\beta^{-1}(t) = \frac{C}{1+\log(t)}$ and we initialize $\nump^{(k)}$ and $\numx^{(k)}$ with $\0\in\mathbb{R}^{Nd}$ and  the Halton point set of appropriate size and dimension, respectively. The other parameters in \cref{alg:annealing}  and the value of the constant $C>0$ in the cooling schedule are tuned for each optimization problem. We ensure that the system has settled in a local minimum at the end of the annealing procedure by monitoring the potential energy trajectory $\left(U({\numx}^{(k)})\right)_{1\leq k \leq T}$. 




\begin{algorithm}
\caption{Underdamped Langevin simulated annealing}\label{alg:annealing}
\begin{algorithmic}[1]
 \STATE {INPUT: $U, \numx^{0},\nump^0, \deltat,T, \beta^{-1}, \gamma$ }
\STATE{$ \alpha := \exp(-\deltat\gamma)$}
\FOR{k = 0 \TO T-1} 	
	\STATE{$\nump^{(k+1/2)}  \gets \nump^{(k)} - \frac{\deltat}{2} \nabla U(\numx^{(k)})$}
	\STATE {$\numx^{(k+1/2)}  \gets \numx^{(k)} + \frac{\deltat}{2}\p^{(k+1/2)} $} 
	\STATE{$\nump^{\prime(k+1/2)} \gets    \alpha \nump^{k+1/2}  +  (1-\alpha^{2})^{1/2} \Rand_{k}, ~ \Rand_{k} \sim \mathcal{N}(\0, \beta^{-1}(\deltat k)\I_{m})$}
	\STATE {$\numx^{(k+1)}  \gets \numx^{(k+1/2)} + \frac{\deltat}{2}\p^{\prime}_{(k+1/2)} $} 
	\STATE{$\nump^{(k+1)}  \gets \p^{\prime(k+1/2)} - \frac{\deltat}{2} \nabla U(\numx^{(k+1)})$}
\ENDFOR
\RETURN{$\arg \min_{\x \in \{ \numx^{(k)}, 1 \leq k \leq T\}} U(\x)$}
\end{algorithmic} \end{algorithm}
\subsection{Simulated annealing scheme for smooth hypersurfaces }\label{ap:sec:annealing:sphere}
For optimization on the sphere, the ``dented'' sphere, and the hyperboloid, we use a constrained version of the Langevin diffusion process \eqref{eq:ULD}, i.e.,  we consider \eqref{eq:ULD} subject to
\begin{equation}\label{eq:constraint}
g(\x) = {\bm 0} \in \mathbb{R}^{N}, 
\end{equation}
and 
\[
\nabla_{\x_{i}} g(\x) \cdot \p_{i} = 0, ~1\leq i \leq N,
\]
where $g = (g_{1},\dots,g_{N})$ is chosen such that the constraint \eqref{eq:constraint} ensures that  $\x_{i}, ~ 1 \leq i \leq N$ are elements of the hypersurface, e.g. 
\[
g_{i}(\x) = \norm{\x_{i}}_{2}^{2} - 1  = 0, ~~1 \leq i \leq N,
\]
in the case of $M=\mathbb{S}^{d-1}$. We use the geodesic Langevin Integrator ``g-BAOAB'' (see \cite{leimkuhler3} for details) in order to numerically integrate the constrained dynamics. The resulting annealing scheme resembles \cref{alg:annealing}. We use a cooling scheme of the same type as in the unconstrained case. In all examples we set the initial velocity of each particle to zero, i.e., $\p^{(0)}_{i} = {\bm 0}, ~~1 \leq i \leq N$. For $M=\mathbb{S}^{2}$ we initialise $\x^{(k)}$ by spherical Fibonacci point sets. In the case of the ``dented'' sphere example, we initialize $\x^{(k)}$ by mapping the spherical Fibonacci point set (for $N=89$) onto the ``dented'' sphere using the projection map
\[
(x_{1},x_{2},x_{3}) \mapsto (x_{1}, {\rm sign}(x_{2}) \sqrt{ (\alpha + x_{1}^2) x_{2}^2}, x_{3}).
\]
In the example of the Poincar\'e disk model we initialize the particles by first generating uniformly distributed points on the disk $\{ (x_{1},x_{2}) : x_{1}^{2}+x_{2}^{2} \leq 4/5 \}$, and then project these points onto the upper sheet of the hyperboloid model using the appropriately defined inverse of the projection \eqref{eq:projection}. During simulation time the additional constraint \eqref{eq:add:constraint} for the compactified hyperboloid is ensured to be (approximately) satisfied by adding the additional energy term $\tilde{U}(\x) = \sum_{i=1}^{N} \tilde{U}_{i}(\x_{i})$, to the energy functional $E_{\rm Gaussian}$, where
\[
 \tilde{U}_{i}( (\x_{i,1},\x_{i,2},\x_{i,3})) = \begin{cases} 0, &\text{if } \abs{\x_{i,3}} \leq  c,\\
 \kappa ( \x_{i,3} - c)^{\alpha}, & \text{otherwise}
 \end{cases}
\]
with $c = \frac{ 1+ r^{2}}{1 - r^{2}}, r =4/5$, and sufficiently large $\alpha>1, \kappa>0$.
As for the torus examples, the values of the remaining parameters in the annealing scheme are chosen problem dependently.


\end{document}